\documentclass[preprint,10pt]{elsarticle}
\usepackage{amssymb}
\usepackage{}
\usepackage{color}
\usepackage{amsfonts}
\usepackage{mathrsfs}
\usepackage{amsthm}
\usepackage{extarrows}
\usepackage{amssymb, amsmath, amscd, latexsym}
\usepackage{booktabs}
\usepackage[colorlinks,linkcolor=blue,anchorcolor=blue,citecolor=blue]{hyperref} 
\usepackage{titlesec}
\usepackage{amsmath} 
\usepackage{algorithm}
\usepackage{algpseudocode}
\usepackage{graphicx}

\usepackage{subcaption} 

\usepackage{titlesec}

\usepackage{multirow}
\usepackage{booktabs}
\usepackage{bm}

\titleformat{\subsection}
  {\normalfont\bfseries} 
  {\thesubsection} 
  {1em} 
  {} 

\large\normalsize
\newtheorem{theorem}{Theorem}[section]

\newtheorem{example}[theorem]{Example}
\newtheorem{remark}[theorem]{Remark}

\numberwithin{equation}{section}

\begin{document}
\newcommand{\Q}{\mathbb{Q}}
\newcommand{\N}{\mathbb{N}}
\newcommand{\Z}{\mathbb{Z}}
\newcommand{\R}{\mathbb{R}}
\newcommand{\D}{\mathbb{D}}
\newcommand{\I}{\mathbb{I}}
\newcommand{\C}{\mathbb{C}}
\newcommand{\T}{\textbf{T}}
\newcommand{\tx}[1]{\quad\mbox{#1}\quad}
\biboptions{sort&compress}%
\definecolor{mgq}{rgb}{0,0,0}
\newcommand{\mgq}{\color{mgq}}
\definecolor{mgq2}{rgb}{0,0,0}
\newcommand{\ma}{\color{mgq2}}
\definecolor{mb}{rgb}{0,0,0}
\newcommand{\mb}{\color{mb}}

\definecolor{yx}{rgb}{0,0,0}
\newcommand{\yx}{\color{yx}}

\begin{frontmatter}
\title{{\bf{ {\color{yx}DDC-PINNs: A Predictor-Corrector Approach Based on Neural Network-Driven Domain Decomposition and Classical ODE Solvers for Time-Dependent PDEs}}}}

%



\author{Xun Yang}
\ead{sicnuyangxun@163.com}

\author{Guanqiu Ma\corref{cor1}}
\ead{gqma@sicnu.edu.cn}

\author{Maohua Ran}
\ead{maohuaran@163.com}

\cortext[cor1]{Corresponding author.}

\address{School of Mathematical Sciences, Sichuan Normal University, Chengdu 610068, China}



\begin{abstract}

When solving time-dependent partial differential equations (PDEs), traditional physics-informed neural networks (PINNs) may encounter several challenges. In particular, standard PINNs do not explicitly account for the temporal evolution order of time-dependent problems during training, which may affect the quality of temporal evolution in time-dependent PDEs. In addition, a single neural network may face difficulties in simultaneously representing different physical behaviors across multiple regions of the computational domain.
To address these issues, we propose a domain-decomposition-based causal PINNs (DDC-PINNs) framework. The term causal refers to the fact that the temporal evolution is performed sequentially through classical ordinary differential equation (ODE) integration, thereby respecting the natural temporal ordering of time-dependent PDEs. The proposed framework enhances spatial approximation through domain decomposition and employs a sequential temporal-evolution strategy for time-dependent problems.
Within this framework, an approximate solution is first obtained using domain-decomposition PINNs. Subsequently, the time-derivative term in the original PDE is retained, while the remaining solution-dependent terms are replaced by the obtained approximation, thereby transforming the original PDE into an auxiliary ODE system. Classical numerical methods for ODEs are then employed to perform temporal evolution without repeated neural-network optimization. As a result, DDC-PINNs decouples spatial approximation from temporal evolution while preserving the temporal evolution order through sequential ODE integration.
Numerical experiments on several benchmark problems demonstrate the effectiveness of the proposed framework and provide proof-of-concept validation of the DDC-PINNs methodology.

\end{abstract}


\end{frontmatter}

\setlength\abovedisplayskip{0.2pt}
\setlength\belowdisplayskip{0.2pt}
\section{Introduction}\label{sec1}
The advent of deep learning has revolutionized computational
methods {\color{mgq} in various domains}, including speech processing,
natural language understanding, and computer vision
\cite{JGu,LECUN,KRIZHEVSKY,HINTON}.
Within this technological landscape,
Raissi et al. proposed PINNs \cite{Raissi}, initiating a new paradigm for applying deep learning to solve forward and inverse problems in PDEs, and providing a potent alternative to conventional numerical PDE solvers. {\color{mgq}Integrating} PDE constraints into the loss function via automatic differentiation, PINNs solve PDEs. {\color{mgq} The advantages of PINNs are mesh-free and structural simplicity.}
{\color{mgq} Thus, PINNs can be used in} various complex PDEs, such as fractional-order \cite{GPang, LingGuo, Zhao2025},  integro-differential\cite{LLu,LeiYuan}, and stochastic PDEs \cite{DZhang,DZhang1}.
Recent research has not only improved the accuracy of PINNs
in solving PDEs \cite{Jeremy,Sifan,Zixue,Sokratis,SWang,LeviD,Yanjie,YiqiGu,Chenxi}, but has also extended their {\color{mgq}applications} to other disciplines \cite{CAIS,MENGX,YChen}.

Despite its many advantages in solving PDEs, PINNs also have certain limitations. First, standard PINNs do not explicitly account for the temporal evolution order of time-dependent problems during training, which may affect long-time prediction accuracy. Second, the training process may become computationally expensive for complex or high-dimensional problems. Moreover, a single neural network may face difficulties in simultaneously representing different physical behaviors across multiple regions of the computational domain. When the solution exhibits strong spatial or temporal variations, a global approximation may become less effective, which can potentially affect prediction accuracy and training robustness.

To address these challenges, researchers have proposed various methods. For instance, Du et al. introduced the Evolutionary Deep Neural Network (EDNN), a novel framework for solving PDEs \cite{YDu}. This algorithm treats the network parameters as variables that evolve over time and updates them dynamically, while embedding constraints such as boundary conditions directly into the network architecture. This design enables EDNN to accurately capture the causality in the solutions of PDEs.
Bruna et al. combined time-parametrized network weights with the
Dirac-Frenkel
variational principle. 
Based on the strong form of the equation, they transformed the PDE into an ODE system, thereby enabling its solution via classical numerical integration methods,
which ensures that the solution possesses causal properties \cite{JBruna,JBerman,YuxiaoWen,Zhang,HuZiqing,Schwerdtner}.
Building on this, Yang et al. first transformed the strong form of the equation into a weak form and then combined it with the Dirac-Frenkel variational principle; this method also captured the causality of the solutions \cite{XYang}; the accuracy of the neural Galerkin method and the weak neural Galerkin method varies depending on the specific case, but their overall accuracy is comparable.
On the other hand, Jung et al. proposed Causal Evolutionary Reinforcement Networks (CEENs)\cite{JungJ}. In this method,
{\color{mgq}time integration is applied to}  both sides of the time-dependent PDEs. By ingeniously employing the trapezoidal rule for numerical time integration, they constructed an equation that {\color{mgq}describes} the causal relationship of the solution. Subsequently, they performed training and optimization at each time step, ultimately obtaining the causality of the solution to this equation.

Jagtap et al. proposed conservative PINNs (cPINNs), which enforce interface continuity
conditions across partitioned domains \cite{Jagtap}, while XPINNs
extend domain decomposition applicability to broader
range of PDEs \cite{ADJagtap}.
Recently, Dolean et al. presented a multi-level domain-decomposition-based PINNs {\mb framework} \cite{DoleanV}. Compared to FBPINNs \cite{Moseley}, it uses a multi-level solution representation. {\color{yx}Different from XPINNs}, {\color{mgq} subdomains of FBPINNs can be overlapped, that avoids the introduction of } extra interface loss terms. This method {\color{mgq} is mainly used in} high-frequency problems.
Kim et al. introduced Initialization-enhanced PINNs (IDPINNs), {\color{mgq} which optimizes} interface continuity in XPINNs and {\color{mgq} contains a strategy for initializing training data}  \cite{Chenhao}.
Hu et al. {\color{yx}presented a non-overlapping Schwarz-type domain decomposition method for physics and equality constrained artificial neural networks} \cite{QifengHu}.
Based on domain decomposition, 
it transforms PDE problems into optimization problems into optimization ones with equality constraints. Its loss function,  built {\color{mgq} by} interface information, {\color{mgq} could be solved by } the adaptive augmented Lagrangian method.
Shang et al. proposed a method that integrates subdomain-decomposition-based stochastic neural networks with overlapping Schwarz preconditioners for PDEs solving \cite{ShangY}. This method combines subdomain-decomposition {\color{mgq} which based on} stochastic neural networks with overlapping Schwarz preconditioners.
Subsequently, more researchers have explored PINNs methods based on domain decomposition techniques.  These developments collectively enhance computational efficiency and theoretical {\color{mgq} basis} for complex system modeling \cite{Ye,ZHu,Hadden}.

To enhance the flexibility of spatial domain processing and improve the treatment of temporal evolution in time-dependent PDEs, we introduce DDC-PINNs, a hybrid framework that combines domain-decomposition PINNs with classical ODE solvers. The method proceeds through two stages.

\begin{enumerate}

\item[(1)]   In the first stage, DDPINNs with enhanced interface treatment are employed to construct an approximation of the solution over decomposed spatial subdomains.

\item[(2)]   In the second stage, the time-derivative term in the PDE is retained, while the remaining solution-dependent terms are replaced by the DDPINNs approximation. This transforms the original PDE into an auxiliary ODE system, which is subsequently evolved using classical numerical integration methods.
\end{enumerate}

The resulting framework separates spatial approximation from temporal evolution. Unlike CEENs~\cite{JungJ}, which require repeated neural-network optimization during time evolution, Neural Galerkin methods~\cite{JBruna,JBerman,YuxiaoWen,Zhang,HuZiqing,Schwerdtner} evolve the neural-network parameters through an ODE system derived from the Dirac--Frenkel variational principle. In contrast, DDC-PINNs first constructs a spatial approximation using DDPINNs and subsequently performs temporal evolution on an auxiliary ODE system obtained from the frozen spatial operator approximation. Consequently, the neural-network parameters remain fixed during the temporal-evolution stage, and no further neural-network optimization or parameter evolution is required.

To clarify the contribution of the proposed framework, we emphasize that DDC-PINNs is not intended as a new time-integration scheme. Neither the classical Runge--Kutta method nor the operator-freezing strategy is individually novel. Instead, the contribution lies in combining

\begin{enumerate}
\item[(1)] domain-decomposition PINNs with enhanced interface regularization;

\item[(2)]  an operator-freezing reformulation that converts the original PDE into an auxiliary ODE system;

\item[(3)]  temporal evolution without repeated neural-network optimization.
\end{enumerate}

This paper is organized as follows: In section \ref{sec2}, PINNs for solving PDEs are reviewed. In section \ref{sec3}, the novel framework DDC-PINNs is introduced. In Section \ref{sec4}, an error analysis of DDC-PINN is presented. In Section \ref{sec5}, the effectiveness of our method is demonstrated through a series of numerical experiments; and our conclusions are presented in Section \ref{sec6}.

{\color{mgq}
\section{Review of PINNs}\label{sec2}
 We review the application of PINNs in solving PDEs in this section.}
We begin by formulating a general system of time-dependent PDEs
\begin{equation}\label{2-1}
u_{t}(\bm{x}, t)+\mathcal{N}[u,\bm{x},t] = 0, \quad \bm{x} \in \Omega, \, t \in [0, T],
\end{equation}
governed by {\color{mgq}initial conditions}
\begin{equation}\label{2-2}
u(\bm{x}, 0) = h(\bm{x}), \quad \bm{x} \in \Omega,
\end{equation}
and {\color{mgq}boundary conditions}
\begin{equation}\label{2-3}
\mathcal{B}[u(\bm{x}, t)] = b(\bm{x}, t), \quad \bm{x} \in \partial\Omega, \, t \in [0, T]{\color{mgq}.}
\end{equation}
Here,  $\Omega $ is a bounded region in $\mathbb{R}^d$($d$ is the dimension of the spatial region), $\mathcal{N}[u,\bm{x},t]$ denotes a linear or nonlinear partial differential operator.
The operator $\mathcal{B}[u(\bm{x}, t)]$ {\color{mgq}represents} the behavior of the solution at the boundary $\partial\Omega \times [0, T]$, and $h(\bm{x})$ represents the initial condition at $t = 0$.

In the PINNs framework, we employ fully connected neural networks to estimate solutions for PDEs.
The {\color{mgq}output of the neural network} $\hat{u}(\bm{x}, t, \theta)$ serves as an approximate solution,
where $\theta$ denotes the {\color{mgq}parameters of the network}. By optimizing $\theta$,
{\color{mgq} it aims} to improve the accuracy of $\hat{u}$ in approximating the true solution.
This optimization is formulated as:
\begin{equation}\label{2-4}
u(\bm{x}, t) \approx \arg\min_{\hat{u} \in F} J(\theta),
\end{equation}
where
\begin{equation}\label{2-5}
F := \{\hat{u}(\cdot, \cdot, \theta) : \theta \in \mathbb{R}^P\} \subset C(\bar{\Omega} \times [0, T]) {\color{mgq}.}
\end{equation}
{\color{mgq} Here, }$P$ denotes the number of weights and bias terms of the neural networks,
and the loss function $J(\theta)$ comprises three key components:
\begin{equation}\label{2-6}
J(\theta) = \lambda _1 \mathcal{L}_r(\theta) + \lambda_2 \mathcal{L}_b(\theta)
 + \lambda_3 \mathcal{L}_0(\theta),
\end{equation}
where $\lambda_1$, $\lambda_2$ and $\lambda_3$ are weighting factors that
 balance the contributions of the PDE residual,
boundary conditions, and initial conditions, respectively.
The components of the loss function of PINNs are defined as follows:
\begin{align}
\mathcal{L}_r(\theta) &= \frac{1}{N_r} \sum_{j=1}^{N_r}
\Big| \hat{u}_{t}(\bm{x}_r^j, t_r^j, \theta)+\mathcal{N}[\hat{u},\bm{x}_r^j, t_r^j]\ \Big|^2, \label{2-7} \\
\mathcal{L}_b(\theta) &= \frac{1}{N_b} \sum_{j=1}^{N_b} \Big| \mathcal{B}[\hat{u}]-(\bm{x}_b^j, t_b^j)
 b(\bm{x}_b^j, t_b^j) \Big|^2,  \label{2-8} \\
\mathcal{L}_0(\theta) &= \frac{1}{N_0} \sum_{j=1}^{N_0} \Big| \hat{u}(\bm{x}_0^j, 0,\theta) -
h(\bm{x}_0^j) \Big|^2, \label{2-9}
\end{align}
where,  $(\bm{x}_r^j, t_r^j) \in \Omega \times [0, T]$, $ (\bm{x}_b^j, t_b^j)
\in \partial\Omega \times [0, T]$, $\bm{x}_0^j \in \Omega $.   $N_r$, $N_b$
and $N_0$ denote the numbers of collocation points for the PDEs residual, boundary conditions and initial conditions, respectively.

{\color{mgq}
\section{DDC-PINNs for solving time-dependent PDEs}\label{sec3}
In this section, we present DDC-PINNs, a novel framework for solving temporal PDEs. First, we introduce DDPINNs as its foundational component.}

\subsection{DDPINNs for solving PDEs}
Domain decomposition is based on the principle of dividing a large domain into smaller subdomains, with each subdomain trained by a separate neural network.
{\color{mgq} Aiming to improve
the accuracy and smoothness of the XPINNs solution at the interface}, we add residual, gradient residual, and continuity terms for both residual gradients and approximate solution gradients at the interface to the total loss function.

Let $\Omega_{1}\cup \cdots\cup\Omega_{i}\cup\cdots \cup \Omega_{j} \cup \cdots  \cup \Omega_{N}= \Omega$ $(1\leq i <j \leq N)$ {\color{mgq}.
If} $\Omega_{i}$, $\Omega_{j}$ are adjacent, then $ \Omega_{i}\cap \Omega_{j}= \partial \Omega_{ij}$; otherwise $ \Omega_{i}\cap \Omega_{j}= \emptyset$.
The training points used for residuals in subdomain $\Omega_{i}$ are denoted {\color{mgq} by} $\{(\bm{x}_q^i, t_q^i)\}_{q=1}^{N_{r_i}}$. ${N_{r_i}}$ denotes the total number of training points on  subdomain $\Omega_{i}$, where $(\bm{x}_q^i, t_q^i)$ $\in \Omega_i \times (0, T)$.
If domain  $\Omega_{i}$ and domain $\Omega_{j}$  are adjacent, $ \{(\bm{x}_q^{ij}, t_q^{ij})\}_{q=1}^{N_{IF_{ij}}} \in \partial \Omega_{ij}$.
$N_{IF_{ij}}$ denotes the number of training points at interface $\partial \Omega_{ij}$. The {\color{mgq} loss function of XPINNs} is defined in the work of Jagtap et al \cite{ADJagtap}. 
In this work, simple non-overlapping domain decompositions are adopted for the benchmark problems in order to illustrate the core idea of the proposed DDC-PINNs framework. The objective of the present study is not to optimize the decomposition topology, but to demonstrate how domain-decomposition-based neural approximations can be coupled with classical ODE solvers for temporal evolution. The solution and residual matching conditions imposed on the interfaces are assumed to be sufficiently smooth so that the corresponding interface loss terms are well defined.
As follows:
\begin{align}
&\mathcal{L}_{IF}(\theta) =\sum_{ij} \frac{1}{N_{IF_{ij}}} \sum_{q=1}^{N_{IF_{ij}}} \Big|
\hat{u}^{i}_{t}(\bm{x}_{q}^{ij}, t_{q}^{ij},\theta)+\mathcal{N}[\hat{u}^{i},\bm{x}_{q}^{ij},
t_{q}^{ij}]-\hat{u}^{j}_{t}(\bm{x}_{q}^{ij}, t_{q}^{ij},\theta)-\mathcal{N}[\hat{u}^{j},\bm{x}_{q}^{ij},t_{q}^{ij}] \Big|^2,\label{2-13}\\
&\mathcal{L}_I(\theta) = \sum_{ij} \frac{1}{N_{IF_{ij}}} \sum_{q=1}^{N_{IF_{ij}}}
\Big| \hat{u}^{i}(\bm{x}_{q}^{ij}, t_{q}^{ij},\theta) -\frac{1}{2}[\hat{u}^{i}(\bm{x}_{q}^{ij},
t_{q}^{ij},\theta)+\hat{u}^{j}(\bm{x}_{q}^{ij}, t_{q}^{ij},\theta)]  \Big|^2, \label{2-14}\\
&J(\theta)_{XP} = \lambda _1 \mathcal{L}_r(\theta) + \lambda_2 \mathcal{L}_b(\theta) +
\lambda_3 \mathcal{L}_0(\theta)+\lambda_{4} \mathcal{L}_{IF}(\theta)+\lambda_{5}\mathcal{L}_I(\theta),\label{2-15}
\end{align}
where $\mathcal{L}_{IF}(\theta)$ {\color{mgq} and} $\mathcal{L}_I(\theta)$ are used to ensure continuity at the interface.

However, this treatment of the interface part is rough and cannot guarantee the smoothness of the interface part \cite{Chenhao}.
To {\color{mgq} deal with the interface part}, we propose DDPINNs based on IDPINNs{\color{mgq}. The }loss function of DDPINNs is constructed as follows:
\begin{equation}\label{2-16}
J(\theta)_{XP_{total}}=J(\theta)_{XP}+ \lambda_6 \mathcal{L}_\nabla \mathfrak{F} (\theta)+
\lambda _7 \mathcal{L}_\nabla u_{c}(\theta) + \lambda_8 \mathcal{L}_\nabla \mathfrak{F}_{c} (\theta)+
\lambda_9 \mathcal{L} \mathfrak{F} (\theta),
\end{equation}
where
\begin{align}
&\mathcal{L}_\nabla \mathfrak{F} (\theta) =\sum_{ij} \sum_{k=i,j}\frac{1}{N_{IF_{ij}}}
\sum_{q=1}^{N_{IF_{ij}}} \Big| \nabla
(\hat{u}^{k}_{t}(\bm{x}_{q}^{ij}, t_{q}^{ij},\theta)+\mathcal{N}[\hat{u}^{k},\bm{x}_{q}^{ij},t_{q}^{ij}])
\Big|^2,  \label{2-17} \\
&\mathcal{L}_\nabla u_{c}(\theta) = \sum_{ij} \frac{1}{N_{IF_{ij}}} \sum_{q=1}^{N_{IF_{ij}}}
 \Big| \nabla \hat{u}^{i}(\bm{x}_{q}^{ij}, t_{q}^{ij},\theta) -\nabla \hat{u}^{j}(\bm{x}_{q}^{ij},
 t_{q}^{ij},\theta)  \Big|^2, \label{2-18} \\
&\mathcal{L}_\nabla \mathfrak{F}_{c} (\theta) =\sum_{ij} \frac{1}{N_{IF_{ij}}} \sum_{q=1}^{N_{IF_{ij}}}
 \Big| \nabla
(\hat{u}^{i}_{t}(\bm{x}_{q}^{ij}, t_{q}^{ij},\theta)+\mathcal{N}[\hat{u}^{i},\bm{x}_{q}^{ij},t_{q}^{ij}])
\label{2-19}\\
&~~~~~~~~~~~~~~~-\nabla(\hat{u}^{j}_{t}(\bm{x}_{q}^{ij}, t_{q}^{ij},\theta)+\mathcal{N}[\hat{u}^{j},\bm{x}_{q}^{ij},t_{q}^{ij}])
  \Big|^2, \nonumber\\
&\mathcal{L}_\mathfrak{F} (\theta) = \sum_{ij} \sum_{k=i,j}\frac{1}{N_{IF_{ij}}} \sum_{q=1}^{N_{IF_{ij}}} \Big|
\hat{u}^{k}_{t}(\bm{x}_{q}^{ij}, t_{q}^{ij},\theta)+\mathcal{N}[\hat{u}^{k},\bm{x}_{q}^{ij},t_{q}^{ij}]
\Big|^2. \label{2-20}
\end{align}

\subsection{DDC-PINNs solving time-dependent PDEs}
Conventional PINNs optimize the solution over the entire space-time domain simultaneously, which may weaken the explicit temporal evolution order in time-dependent problems. To address this issue, we propose the DDC-PINNs method, which is built on a domain-decomposition framework. This method retains the flexibility of DDPINNs for spatial approximation while introducing a sequential temporal-evolution procedure through an auxiliary ODE system.

{\color{mgq}Using} the approximate solution $\hat{u}$ obtained by DDPINNs, we establish the approximation:
\begin{equation}\label{2-21}
\mathcal{N}[u,\bm{x},t] \approx \mathcal{N}[\hat{u},\bm{x},t].
\end{equation}
{\color{mgq} By this approximation, the original problem \eqref{2-1} -\eqref{2-3} 
transforms into the following initial-value problem of an ODE system
\begin{align}
&u_t(\bm{x}, t) \approx -\mathcal{N}[\hat{u},\bm{x},t],~~ \bm{x} \in \Omega, t \in [0, T], \label{2-22} \\
&u(\bm{x}, 0) = h(\bm{x}),~~\bm{x} \in \Omega. \label{2-23}
\end{align}
{\color{mgq} The system} \eqref{2-22}-\eqref{2-23} can be solved {\color{mgq} by} classical numerical methods such as the Euler scheme or Runge-Kutta methods. The resulting two-stage framework is referred to as DDC-PINNs.

\begin{remark}
The purpose of DDC-PINNs is not to introduce a new ODE integrator.
Instead, the framework combines domain-decomposition PINNs and classical time integration in a two-stage manner.
The DDPINNs approximation is first used to construct a frozen approximation of the nonlinear spatial operator. The resulting auxiliary ODEs system is then evolved independently of the neural-network parameters.
This is different from Neural Galerkin methods, where the network parameters themselves evolve according to an ODE system, and from CEENs, where neural-network optimization is repeatedly performed during temporal evolution.
Consequently, DDC-PINNs separates spatial approximation from temporal evolution and avoids repeated neural-network retraining in the second stage.
\end{remark}

Figure \ref{fig:1} illustrates the DDC-PINNs framework for solving time-dependent PDEs, where  $\hat{u}$ denotes the approximate solution obtained by DDPINNs and  $\tilde{u}$ denotes the solution obtained by classical numerical integration. The framework reflects the sequential temporal evolution of the solution. For clarity, the main implementation procedure corresponding to Figure \ref{fig:1} is summarized in Algorithm \ref{alg1}.

\begin{figure}[h]
    \centering
    \includegraphics[scale=0.8]{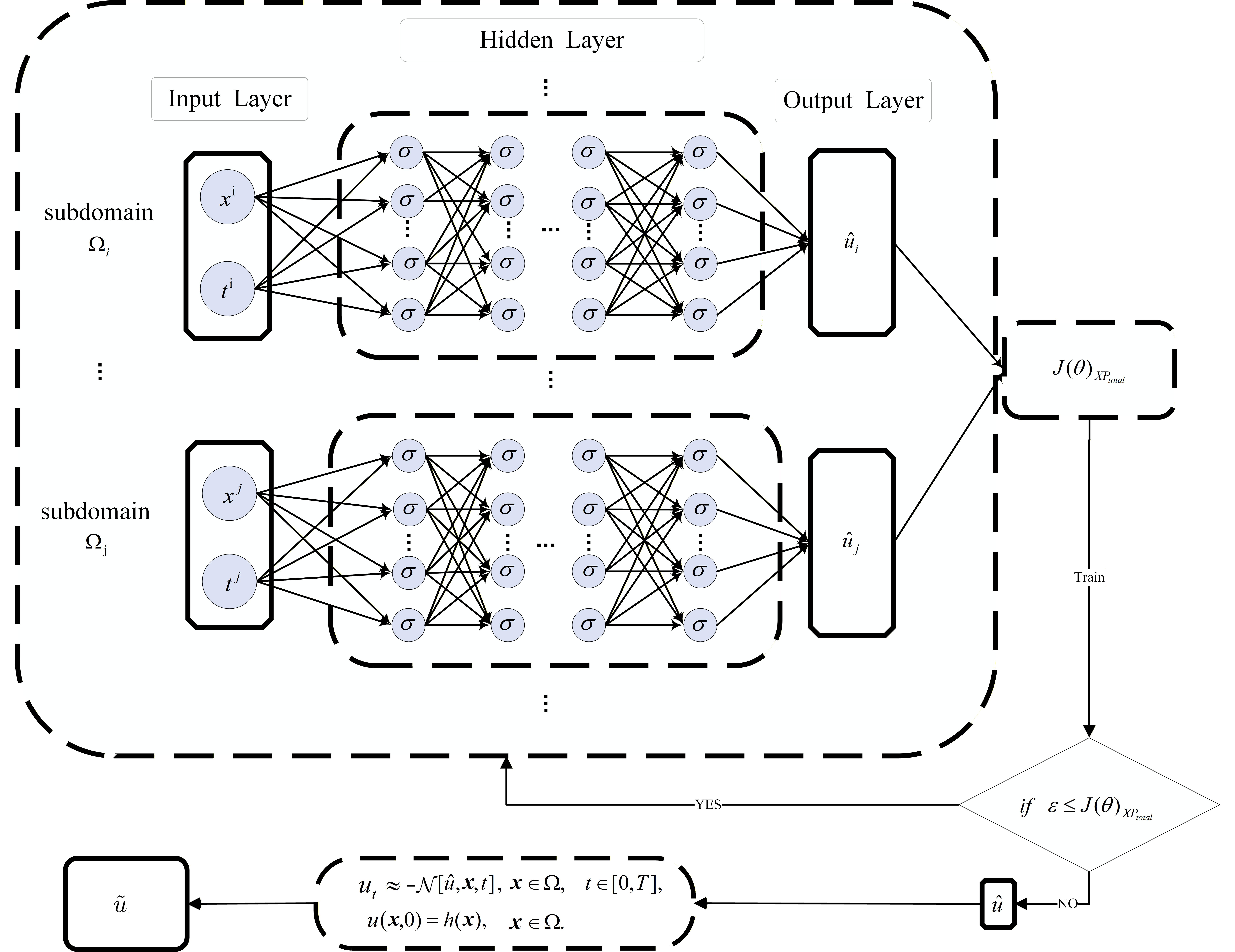}
    \caption{Schematic illustration of the DDC-PINNs framework. Stage I employs DDPINNs on decomposed subdomains to construct an approximation $\hat{u}$ and the corresponding operator approximation $N[\hat{u},x,t]$. Stage II transforms the original PDE into an auxiliary ODE system and performs temporal evolution using a classical numerical integrator (e.g.,Runge-Kutta) without further neural-network optimization.}
    \label{fig:error1}
      \label{fig:1}
\end{figure}

\begin{algorithm}[h]
\caption{Implementation procedure of DDC-PINNs}\label{alg1}
\begin{algorithmic}[1]
\State Generate collocation sets in each subdomain and on interfaces.
\State Train the DDPINNs by minimizing the total loss $J(\theta)_{XP_{total}}$ defined in Eq~\eqref{2-16}.
\State Construct the domain-decomposition approximation $\hat{u}=\sum_i \hat{u}_i \mathbf{1}_{\Omega_i}$.
\State Compute the frozen operator $\mathcal{N}[\hat{u},\bm{x},t]$.
\State Solve the auxiliary ODE system
$\tilde{u}_t(\bm{x},t)=-\mathcal{N}[\hat{u},\bm{x},t],\quad \tilde{u}(\bm{x},0)=h(\bm{x})$ using a classical numerical method (e.g.,  Runge-Kutta).
\end{algorithmic}
\end{algorithm}

\section{Consistency and Error Analysis}\label{sec4}
In this subsection, we analyze the error associated with the proposed DDC-PINNs framework. This error depends primarily on a pair of key factors: the initial approximation accuracy provided by the DDPINNs, and the numerical error introduced by the subsequent ordinary differential equation integrator. 
Following the error decomposition introduced for XPINNs  \cite{ADJagtap}, we analyze the error sources in the DDPINNs approximation.
\subsection{Error Sources in the DDPINNs Approximation}
Let the computational domain be partitioned into  $N_{sd}$  nonoverlapping subdomains  $ \{ \Omega_{q}\}^{N_{sd}}_{q=1}$, In each subdomain a separate neural network $\hat{u}=\Sigma_{q}\hat{u}_{q} \mathbf{1}_{\Omega_{q}}$ is trained ($\mathbf{1}_{\Omega_{q}}$ is an indicator function). 
The total error between the exact solution $u$ and the DDPINNs approximation can be decomposed as follows.
The total error can be decomposed as
\begin{equation}\label{3-1}
\mathcal{E}_{\mathrm{total}} = E(\mathcal{E}_{\mathrm{app}}, \mathcal{E}_{\mathrm{gen}},\mathcal{E}_{\mathrm{opt}})
\end{equation}
where
\begin{itemize}
    \item \(\mathcal{E}_{\mathrm{app}}\) is the \textbf{approximation error}, arising from the limited capacity of the neural network to represent the exact solution. It can be reduced by increasing the depth/width of the network or by using adaptive activation functions.
    \item \(\mathcal{E}_{\mathrm{gen}}\) is the \textbf{generalization error}, caused by the finite number and distribution of collocation points (residual, boundary, and interface points). This error depends on the sampling strategy and can be mitigated by placing more points in regions with complex solution behaviour or near interfaces.
    \item \(\mathcal{E}_{\mathrm{opt}}\) is the \textbf{optimization error}, originating from the non-convex nature of the loss landscape and the difficulty of reaching the global minimum. It is influenced by the choice of optimizer, learning rate, and the number of training iterations. $E$ is an abstract function that depends on $\mathcal{E}_{\mathrm{app}}$, $\mathcal{E}_{\mathrm{gen}}$ and $\mathcal{E}_{\mathrm{opt}}$.
\end{itemize}

\subsection{Error Propagation in the Reduced ODEs System}
For the theoretical analysis presented below, the exact solution $u$ and the DDPINNs approximation $\hat{u}$ are assumed to possess the regularity required by the subsequent derivations, so that all derivatives appearing in the analysis are well defined.
Let $X$ and $Y$ denote Banach spaces such that
$
u(\cdot,t)\in X
$
and
$
\mathcal{N}[u,\cdot,t]\in Y
$
for each $t\in[0,T]$. The nonlinear spatial operator is regarded as a mapping
$
\mathcal{N}:X\rightarrow Y.
$
For example, for the diffusion problem considered in Example~\ref{ex:1}, one may take
$
X=H^1(\Omega),\quad Y=L^2(\Omega),
$
where $\Omega=[-1,1]$.
Unless otherwise stated, $\|\cdot\|_X$ and $\|\cdot\|_Y$ denote the corresponding norms in these spaces.

Let $u$ be the exact solution of the original PDE
\begin{equation}\label{3-2}
u_t(\bm{x}, t) + \mathcal{N}[u,\bm{x},t] = 0,\qquad u(\bm{x},0)=h(\bm{x}),
\end{equation}
and let $\hat{u}(\bm{x}, t)$ be the DDPINNs approximation. In the DDC-PINNs framework, we solve the initial value problem
\begin{equation}\label{3-3}
\phi_t(\bm{x}, t) = -\mathcal{N}[\hat{u},\bm{x},t],\qquad \phi(\bm{x},0)=h(\bm{x}),
\end{equation}
where $\phi(\bm{x}, t)$ denotes the DDC-PINNs solution. Define the error $e(\bm{x}, t) = \phi(\bm{x}, t)- u(\bm{x}, t)$. Using \eqref{3-3}-\eqref{3-2}, we obtain

\begin{equation}\label{3-4}
e_t (\bm{x}, t)= -\mathcal{N}[\hat{u},\bm{x},t]+ \mathcal{N}[u,\bm{x},t] = \mathcal{N}[u,\bm{x},t]  - \mathcal{N}[\hat{u},\bm{x},t].
\end{equation}
Let $\delta (\bm{x}, t)= u (\bm{x}, t)- \hat{u} (\bm{x}, t)$ be the DDPINNs approximation error. Then $u(\bm{x}, t) = \hat{u}(\bm{x}, t) + \delta (\bm{x}, t)$ and the error equation becomes
\begin{equation}\label{3-5}
e_t(\bm{x}, t) = \mathcal{N}[\hat{u}+\delta,\bm{x},t] - \mathcal{N}[\hat{u},\bm{x},t]. 
\end{equation}
For any time $t$, assume that the exact solution $u(\cdot,t)$ and the DDPINNs approximation $\hat{u}(\cdot,t)$ belong to a bounded convex set $B\subset X$ (e.g., $B={w\in H^1(\Omega): \|w\|_{H^1}\le M}$). On this bounded set, the operator $\mathcal{N}:X\rightarrow Y$ is assumed to be locally Lipschitz continuous. More precisely, there exists a constant $L>0$ (depending on $B$ but independent of $t$) such that for all $t\in[0,T]$ and any two functions $w_1(\cdot,t), w_2(\cdot,t)\in B$, we have
\begin{equation}\label{3-6}
\left\|\mathcal{N}[w_1(\cdot,t),\bm{x},t]-\mathcal{N}[w_2(\cdot,t),\bm{x},t]\right\|_{Y}
\le
L\left\|w_1(\cdot,t)-w_2(\cdot,t)\right\|_{X}.
\end{equation}
For the benchmark problems considered in this work, the exact solution and the DDPINNs approximation are assumed to remain in $B$ for all $t\in[0,T]$. Here the norms $\|\cdot\|_X$ and $\|\cdot\|_{Y}$ depend only on the spatial variable $\bm{x}$, while the time variable $t$ is treated as a parameter.

We integrate the error equation $e_t(\bm{x}, t) = \mathcal{N}[\hat{u}+\delta,\bm{x},t] - \mathcal{N}[\hat{u},\bm{x},t]$ over the interval $[0,s]$ to obtain
\begin{equation}\label{3-7}
e(\bm{x}, s) = \int_0^s \bigl( \mathcal{N}[\hat{u}+\delta,\bm{x},t] - \mathcal{N}[\hat{u},\bm{x},t] \bigr) dt.
\end{equation}
Furthermore, we obtain
\begin{equation}\label{3-8}
\left\|e(\bm{x}, s)\right\|_{Y} =\left\|\int_0^s \mathcal{N}[\hat{u}+\delta,\bm{x},t] - \mathcal{N}[\hat{u},\bm{x},t]dt\right\|_{Y} \le \int_0^s \left\|\mathcal{N}[\hat{u}+\delta,\bm{x},t] - \mathcal{N}[\hat{u},\bm{x},t]\right\|_{Y} dt.
\end{equation}
Applying the Lipschitz condition \eqref{3-5} pointwise for each $s$ yields
\begin{equation}\label{3-9}
\left\|e(\bm{x}, s)\right\|_{Y} \le L \int_0^s \left\|\delta(\bm{x}, t)\right\|_{X} dt.
\end{equation}
Here, $\delta(\bm{x}, t) = u(\bm{x}, t) - \hat{u}(\bm{x}, t)$ represents the approximation error between the DDPINNs approximated solution and the exact solution.

Inequality \eqref{3-9} shows that the error introduced by replacing
$\mathcal{N}[u,\bm{x},t]$ with $\mathcal{N}[\hat{u},\bm{x},t]$
is controlled by the DDPINNs approximation error, thereby providing
a consistency justification for the proposed DDC-PINNs framework.

For nonlinear problems (e.g., Burgers or KdV), $\mathcal{N}$ is locally Lipschitz on bounded sets. This follows from the fact that the nonlinear terms involve products of the solution and its spatial derivatives. When both the exact solution $u$ and the DDPINNs approximation $\hat{u}$ remain bounded in the chosen function space, these nonlinear operators satisfy a local Lipschitz condition on the corresponding bounded subset.
Since the exact solution $u$ is bounded (by well-posedness) and the DDPINNs approximation $\hat{u}$ is accurate, both $u$ and $\hat{u}$ lie in a common bounded set where the Lipschitz condition holds uniformly. Consequently, the error estimate
$\|e(t)\|_Y \le L\int_0^t \|\delta(s)\|_X \, ds$.
ensures that the DDC-PINNs solution depends continuously on $\hat{u}$.
Time integration is performed with the explicit fourth-order Runge–Kutta method ($\Delta t=0.001$). Since the right-hand side of the auxiliary ODE system is independent of $\phi$, the temporal evolution stage does not involve repeated
updates of neural-network parameters. The stability of DDC-PINNs depends on two components:
\begin{itemize}
\item[(i)] the approximation accuracy of DDPINNs;
\item[(ii)] the numerical integration error of the ODE solver.
\end{itemize}
The estimate \eqref{3-9} indicates that the DDC-PINNs solution depends continuously on the DDPINNs approximation. Consequently, improvements in the first-stage approximation directly reduce the temporal evolution error.  A rigorous stability analysis of the coupled DDPINNs--ODE framework is beyond the scope of the present work and will be investigated in future research.

\section{ Numerical experiments}\label{sec5}
In this section, we evaluate the performance of the proposed DDC-PINNs framework on four benchmark problems. Comparisons are conducted with PINNs, XPINNs, IDPINNs, CEENs, and the Neural Galerkin method. The selected baseline methods are chosen according to their relevance to the main components of the proposed DDC-PINNs framework. XPINNs and IDPINNs are included because they are representative domain-decomposition-based PINNs, allowing us to evaluate the effect of the domain-decomposition component. 
CEENs and Neural Galerkin methods are considered because they incorporate temporal-evolution mechanisms and therefore provide meaningful comparisons for assessing the sequential temporal-evolution component of DDC-PINNs.
Although modern operator-learning approaches such as DeepONets, FNOs, and transformer-based PDE solvers have shown strong performance in many PDE problems, they are designed for different problem settings and are not directly focused on domain decomposition or sequential ODE-based temporal evolution.

To ensure fairness in the comparisons, identical neural network architectures and training settings were used for each benchmark problem whenever possible. The Adam optimizer was used for all experiments. The detailed hyperparameter settings are summarized in Table~\ref{tab:0}.

It should be noted that DDC-PINNs, XPINNs, and IDPINNs are domain-decomposition-based methods and therefore employ multiple subnetworks associated with different spatial subdomains. In contrast, CEENs and Neural Galerkin methods use a single neural network following their original formulations. For each subdomain in DDC-PINNs, the same network architecture as that of the corresponding single-network methods is adopted. Therefore, the use of multiple subnetworks originates from the domain-decomposition strategy rather than from intentionally increasing the capacity of an individual neural network.

In the experiments involving DDC-PINNs and the Neural Galerkin method, the fourth-order Runge–Kutta method is employed with a time-step size of 0.001. For CEENs, the time-step size is likewise set to 0.001, and the error tolerance of the stopping criterion is set to $10^{-5}$.
The loss weights $\lambda$ are chosen empirically to balance the contributions of the PDE residual, boundary-condition, initial-condition, and interface-related losses during training\cite{Rohrhofer2023}, and the corresponding values are reported in each subsection .

\begin{table}[htbp]
\centering
\caption{Neural Network Hyperparameter Settings}
 \label{tab:0}
\begin{tabular}{|c|c|c|c|c|c|c|}
\hline
Example & Methods & Hidden Layers & Neurons & Activation & Learning Rate & Iterations \\
\hline
\ref{ex:1} & All methods & 4 & 50 & $tanh$ & 0.001 & 20000 \\
\hline
\ref{ex:2} & All methods & 4 & 50 & $tanh$ & 0.001 & 20000 \\
\hline
\ref{ex:3}  & All methods & 2 & 20 & $tanh$ & 0.001 & 20000 \\
\hline
\ref{ex:4}  & All methods & 4 & 32 & $tanh$ & 0.001 & 20000 \\
\hline
\end{tabular}
\end{table}

The relative $L^2$  error is defined as
\begin{align}
&\text{$L^2$-error} = \frac{\sqrt{\sum\limits_{j=1}^N |\hat{u}(\bm{x}_j, t_j, \cdot) - u_{\text{ref}}(\bm{x}_j, t_j)|^2}}{\sqrt{\sum\limits_{j=1}^N |u_{\text{ref}}(\bm{x}_j, t_j)|^2}}, \label{4-2} 
\end{align}
here, {\ma $N$, $\hat{u}$ and $u_{\text{ref}}$ denote the total number of computational points, the approximate solution, and the exact solution or reference solution, respectively.}

To evaluate the robustness of the proposed method with respect to random initialization, all experiments were independently repeated five times using random seeds ${0,1,2,3,4}$. The quantitative results reported in all tables are presented as mean $\pm$ standard deviation of the relative $L^2$ errors and runtimes over these independent runs. Since the variations among different runs are relatively small, all solution visualizations and point-wise error plots shown in the figures are generated using the results corresponding to seed $=0$, which is representative of the overall numerical behavior observed across different runs.

All runtime statistics are computed on the same hardware, a consumer-grade laptop (Lenovo XiaoXin Air 15 ALC 2021) equipped with an AMD Ryzen 7 5700U CPU (8 cores, 16 threads), 16 GB of DDR4 RAM (2666 MHz, onboard), and integrated AMD Radeon Graphics (Vega 8).
All reported runtimes correspond to the total wall-clock time required to complete the entire computational workflow of each method under the same hardware environment. This includes all stages necessary to obtain the final numerical solution and corresponding error measurements. The reported runtimes are intended as practical performance indicators under a unified experimental setting.

\subsection{Diffusion problem}

\begin{example}\label{ex:1}
This example considers a one-dimensional diffusion process with a source term.
\begin{align}
    &u_t=u_{xx}
     +e^{-t}(\pi^{2}\sin(\pi x)-\sin(\pi x)) , \quad x \in [-1, 1], t \in [0, 1], \label{4-4}\\
    &u(-1, t) = u(1, t) = 0, \quad t \in [0, 1], \label{4-5}\\
    &u(x, 0) = \sin(\pi x),  \quad x \in [-1, 1].\label{4-6}
\end{align}
The exact solution of problem \eqref{4-4}-\eqref{4-6} is $u(x, t) =e^{-t}\sin(\pi x)$.
\end{example}

In the experimental setup, the spatial domain $[-1, 1]$ is partitioned into two subdomains $[-1, 0]$ and $[0, 1]$, leading to two space-time regions$[-1, 0] \times [0, 1]$ and $[0, 1] \times [0, 1]$.
The loss function weights are assigned as follows:
$\lambda_{1} = 1$, $\lambda_{2} = 22$, $\lambda_{3} = 1$, $\lambda_{4} = 1$, $\lambda_{5} = 10$,$\lambda_{6} = 0.001$, $\lambda_{7} = 10$, $\lambda_{8} = 0.001$, $\lambda_{9} = 1$.
Training Point Allocation for XPINNs, IDPINNs, and DDC-PINNs (First Stage):
\begin{itemize}
\item Each subdomain contains 3000 points.
\item An additional 1000 points are allocated to the subdomain interface.
\item 200 points enforce boundary conditions.
\item 80 points enforce initial conditions.
\end{itemize}
Here, for the CEENs,  neural Galerkin methods, and DDC-PINNs methods (Stages 2), the spatial grid consists of $200$ points, and the time step is $0.001.$

Figure~\ref{fig:2} presents the exact solution, the DDC-PINNs prediction, and the corresponding two-dimensional and three-dimensional visualizations of the point-wise errors for problem \eqref{4-4}-\eqref{4-6}. The solution exhibits the expected diffusion behavior, with the amplitude gradually decreasing over time.
In Figure  \ref{fig:3}, the DDC-PINNs solution (blue solid line) is compared with the exact solution (red dashed line) at different time instances.

In Table \ref{tab:1}, we summarize the relative $L^2$ errors and runtimes obtained by PINNs and the domain-decomposition methods (XPINNs, IDPINNs, and DDC-PINNs) for Example~\ref{ex:1}. Among the compared methods, DDC-PINNs produce the smallest mean relative $L^2$ error. In terms of runtime, DDC-PINNs require more computational time than PINNs, XPINNs, and IDPINNs in this example.
In Table \ref{tab:2}, we compare the relative $L^2$ errors and runtimes of the Neural Galerkin method, the CEENs method, and the DDC-PINNs method. For this example, DDC-PINNs achieve a lower mean relative $L^2$ error than both Neural Galerkin and CEENs. Regarding runtime, the Neural Galerkin method is faster than DDC-PINNs, whereas CEENs require a longer computational time.

\begin{table}[h]
\centering
\caption{Relative $L^2$ error comparison among PINNs, XPINNs, IDPINNs, and DDC-PINNs for Example~\ref{ex:1}.}
 \label{tab:1}
\begin{tabular}{ccccc}
\toprule
\textbf{Methods} & PINNs & XPINNs & IDPINNs  & DDC-PINNs \\
\midrule
Relative \(L^2\) error  &  $(2.82\pm 0.53)\%$ &  \( (8.45\pm 5.81 )\% \) & \((1.18\pm 0.31)\%\) & \((0.78\pm 0.17)\%\) \\
Runtime (s)
&
$1004.46\pm9.55 $
&
$1304.07\pm 27.26$
&
$1823.54\pm 21.32$
&
$1876.42\pm 24.19$
\\
\bottomrule
\end{tabular}
\end{table}

\begin{table}[h]
\centering
\caption{Accuracy and runtime comparison of Neural Galerkin, CEENs, and DDC-PINNs for Example~\ref{ex:1}.}
\label{tab:2}
\begin{tabular}{cccc}
\toprule
Methods
&
Neural Galerkin
&
CEENs
&
DDC-PINNs
\\
\midrule
Relative $L^2$ error
&
$(1.45\pm 1.11)\%$
&
$(15.52\pm 0.28)\%$
&
$(0.78\pm 0.17)\%$
\\
Runtime (s)
&
$1284.02\pm 51.79$
&
$2613.67 \pm 57.94 $
&
$1876.42\pm 24.19$
\\
\bottomrule
\end{tabular}
\end{table}

\begin{figure}[H]
    \centering
    \begin{subfigure}[b]{0.4\textwidth}
        \centering
        \includegraphics[width=\textwidth]{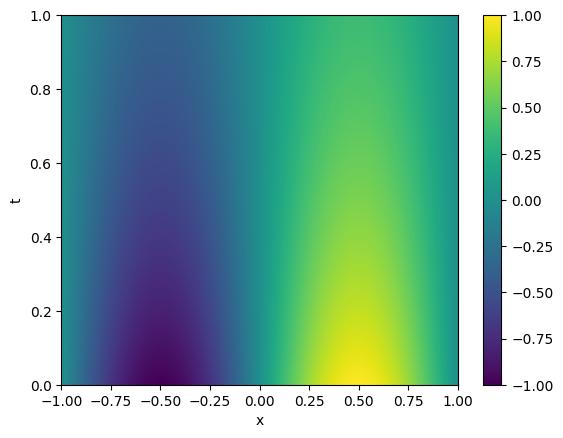}
        \caption{Exact solution}
    \end{subfigure}
    \vspace{0.5cm}
    \begin{subfigure}[b]{0.4\textwidth}
        \centering
        \includegraphics[width=\textwidth]{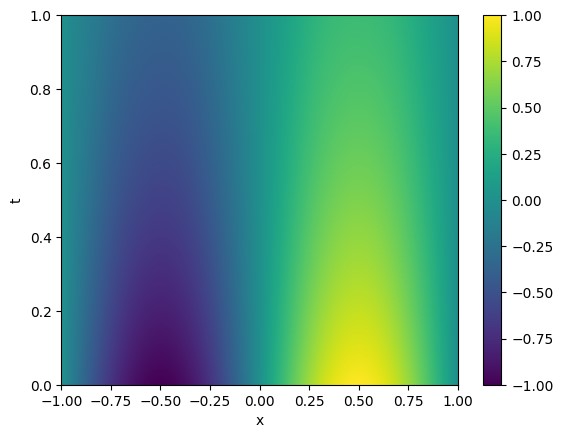}
        \caption{DDC-PINNs solution}
    \end{subfigure}
    \hfill
    \begin{subfigure}[b]{0.4\textwidth}
        \centering
        \includegraphics[width=\textwidth]{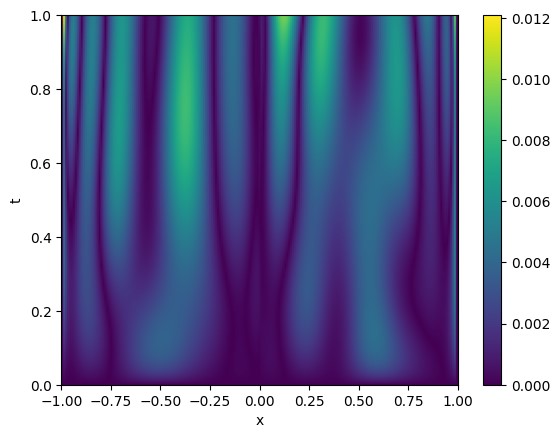}
        \caption{Point-wise errors}
    \end{subfigure}
    \vspace{0.5cm}
    \begin{subfigure}[b]{0.4\textwidth}
        \centering
        \includegraphics[width=\textwidth]{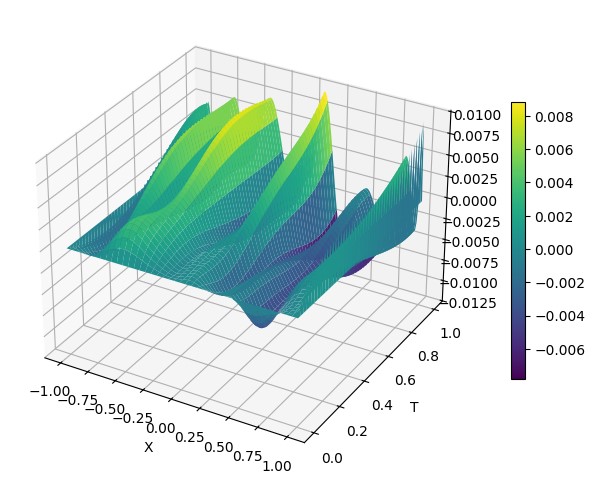}
        \caption{3D point-wise errors}
    \end{subfigure}
     \caption{Comparison of the exact solution, the DDC-PINNs prediction, and the corresponding point-wise errors for Example~\ref{ex:1}. Panels (c) and (d) present the two-dimensional point-wise absolute errors and the three-dimensional point-wise errors, respectively. The solution exhibits the expected diffusion behavior, with the amplitude gradually decreasing over time.}\label{fig:2}
\end{figure}


\begin{figure}[H]
    \centering
    \begin{subfigure}[b]{1\textwidth}
        \centering
        \includegraphics[width=\textwidth]{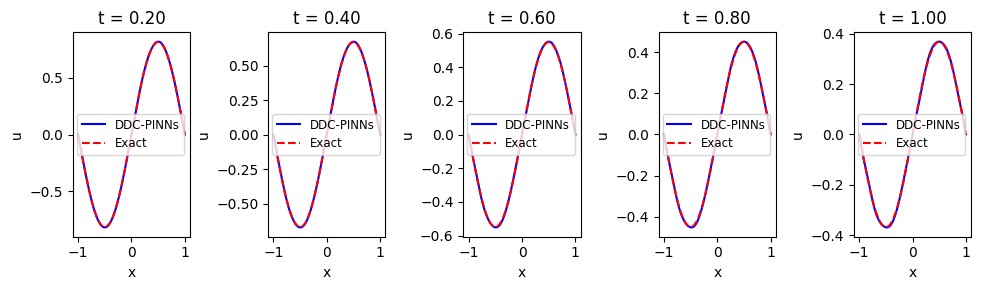}
    \end{subfigure}
    \caption{
    Comparison between the exact solution (red dashed line) and the DDC-PINNs solution (blue solid line) at representative time instants for Example \ref{ex:1}.
The close agreement of the two curves indicates that the proposed framework accurately reproduces the temporal evolution of the diffusion process.}  \label{fig:3}
\end{figure}

\subsection{ Inviscid Burgers' equation}
\begin{example}\label{ex:2}
 This example considers the inviscid Burgers' equation, which is commonly used as a benchmark nonlinear convection–diffusion problem.
\begin{align}
    &u_t =- uu_x,   \quad x \in [-1, 1], ~~t \in [0, 1], \label{4-8}\\
    &u(-1, t) = u(1, t) , \quad t \in [0, 1], \label{4-9}\\
    &u(x, 0) = 0.3\exp(-9x^2)+1,  \quad x \in [-1, 1].\label{4-10}
\end{align}
\end{example}

In this example, the reference solution is computed {\ma by the} finite difference method with central differencing for both temporal ($\Delta t=0.0001$) and spatial ($\Delta x=0.01$) discretizations. 
The spatial domain \([-1,1]\) is divided into two subdomains, $[-1,0]$ and $[0,1]$. 
The weights in the loss function are set as $\lambda_{1}=5$, $\lambda_{2}=20$, $\lambda_{3}=10$, $\lambda_{4}=0.01$, $\lambda_{5}=0.0002$, $\lambda_{6}=0.01$, $\lambda_{7}=0.01$, $\lambda_{8}=0.01$, $\lambda_{9}=0.01$. 
Regarding the training point allocation for XPINNs, IDPINNs, and DDC-PINNs (first stage):
\begin{itemize}
\item Each subdomain is assigned 3000 points.
\item The subdomain interface is given an additional 2000 points.
\item Boundary conditions are enforced by 200 points.
\item Initial conditions are established through 40 points.
\end{itemize}
Here, for the CEENs, neural Galerkin method, and DDC‑PINNs (stages 2), the number of spatial discretization points is 200 and the time step size is 0.001.


Figure~\ref{fig:4} presents the reference solution, the DDC-PINNs prediction, and the corresponding two-dimensional and three-dimensional visualizations of the point-wise errors for problem \eqref{4-8}-\eqref{4-10}. The solution exhibits nonlinear wave propagation driven by the convective term.
Figure \ref{fig:5} presents a comparison between the DDC-PINNs solution (blue solid line) and the reference solution (red dashed line) at different time instances.

In Table \ref{tab:3}, we summarize the relative $L^2$ errors and runtimes obtained by PINNs and the domain-decomposition methods (XPINNs, IDPINNs, and DDC-PINNs) for Example~\ref{ex:2}. Among the compared methods, DDC-PINNs produce the smallest mean relative $L^2$ error, while IDPINNs achieve a comparable level of accuracy. In terms of runtime, DDC-PINNs require slightly more computational time than PINNs, XPINNs, and IDPINNs in this example.
In Table \ref{tab:4}, we compare the relative $L^2$ errors and runtimes of the Neural Galerkin method, the CEENs method, and the DDC-PINNs method. For this example, DDC-PINNs achieve the lowest mean relative $L^2$ error among the three methods. In terms of runtime, DDC-PINNs require less computational time than both the Neural Galerkin method and CEENs.

\begin{table}[h]
\centering
\caption{Relative $L^2$ error comparison among PINNs, XPINNs, IDPINNs, and DDC-PINNs for Example~\ref{ex:2}.}
 \label{tab:3}
\begin{tabular}{ccccc}
\toprule
\textbf{Methods} & PINNs & XPINNs & IDPINNs  & DDC-PINNs \\
\midrule
Relative \(L^2\) error  &   \((1.15\pm 0.46)\%\) & \((3.63\pm 0.98)\%\) & \((0.39\pm 0.09)\%\) &  \((0.38\pm 0.02)\%\)\\
Runtime (s)
&
$505.62\pm 5.28 $
&
$740.42\pm 2.86$
&
$994.82\pm 11.83$
&
$1013.43\pm 11.57$
\\
\bottomrule
\end{tabular}
\end{table}

\begin{table}[h]
\centering
\caption{Accuracy and runtime comparison of Neural Galerkin, CEENs, and DDC-PINNs for Example~\ref{ex:2}.}
\label{tab:4}
\begin{tabular}{cccc}
\toprule
Methods
&
Neural Galerkin
&
CEENs
&
DDC-PINNs
\\
\midrule
Relative $L^2$ error
&
$(0.56\pm 0.11)\%$
&
$(10.26\pm 0.01)\%$
&
$ (0.38\pm 0.02)\%$
\\
Runtime (s)
&
$1105.44\pm 69.02$
&
$1199.36\pm 42.72$
&
$1013.43\pm 11.57$
\\
\bottomrule
\end{tabular}
\end{table}

\begin{figure}[H]
    \centering
    \begin{subfigure}[b]{0.4\textwidth}
        \centering
        \includegraphics[width=\textwidth]{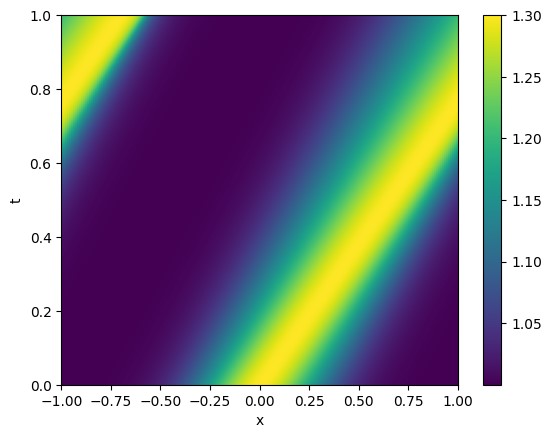}
        \caption{Reference solution}
    \end{subfigure}
    \vspace{0.5cm}
    \begin{subfigure}[b]{0.4\textwidth}
        \centering
        \includegraphics[width=\textwidth]{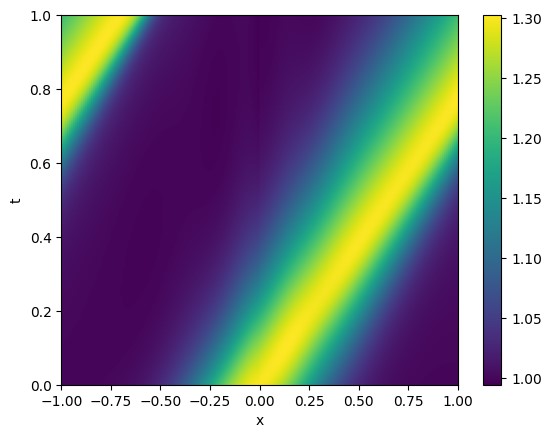}
        \caption{DDC-PINNs solution}
    \end{subfigure}
    \vspace{0.5cm}
    \begin{subfigure}[b]{0.4\textwidth}
        \centering
        \includegraphics[width=\textwidth]{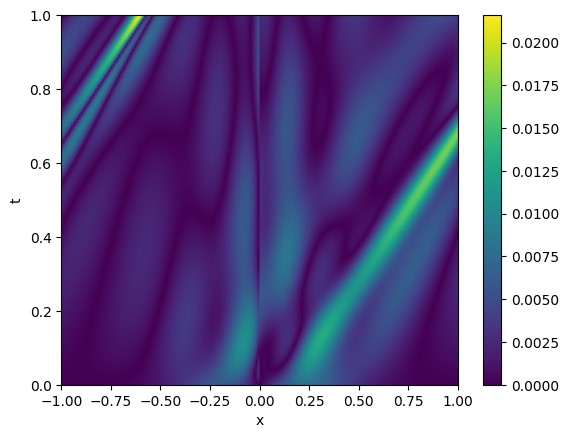}
        \caption{Point-wise errors}
    \end{subfigure}
      \vspace{0.5cm}
    \begin{subfigure}[b]{0.4\textwidth}
        \centering
        \includegraphics[width=\textwidth]{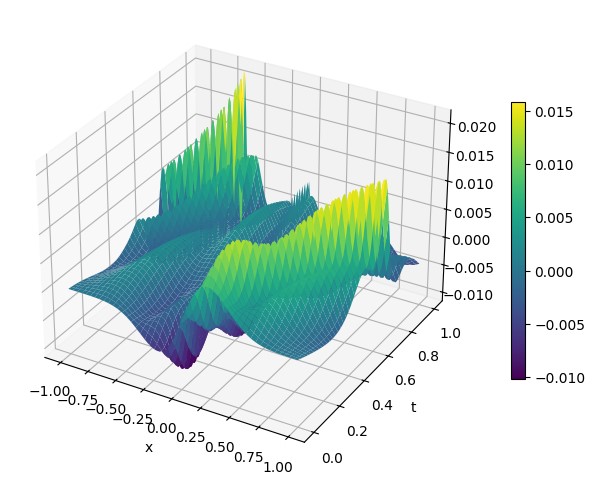}
        \caption{3D point-wise errors}
    \end{subfigure}
    \caption{Comparison of the reference solution, the DDC-PINNs prediction, and the corresponding point-wise errors for Example~\ref{ex:2}. Panels (c) and (d) present the two-dimensional point-wise absolute errors and the three-dimensional point-wise errors, respectively. The solution exhibits nonlinear wave propagation driven by the convective term.} \label{fig:4}
\end{figure}

\begin{figure}[H]
    \centering
    \begin{subfigure}[b]{1\textwidth}
        \centering
        \includegraphics[width=\textwidth]{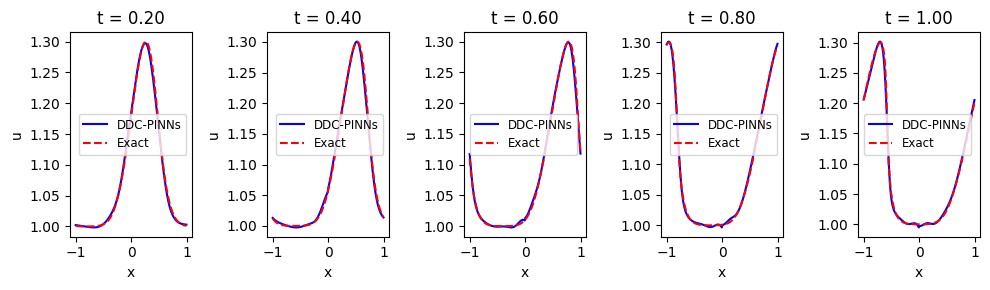}
    \end{subfigure}
    \caption{ Comparison between the reference solution (red dashed line) and the DDC-PINNs solution (blue solid line) at representative time instants for Example \ref{ex:2}.
The results show that the proposed framework successfully tracks the temporal evolution of the nonlinear wave profile.}      \label{fig:5}
\end{figure}

\subsection{Korteweg-de Vries equation(KdV)}
\begin{example}\label{ex:3}
This example considers the KdV equation, which describes nonlinear dispersive wave propagation and solitary-wave dynamics \cite{Alexander}.
\begin{align}
    &u_t+uu_{x}+u_{xxx}=0, \quad (x,t) \in  [-40, 40] \times [0, 1], \label{4-11}\\
    &u(-40,t) =u(40,t),     \quad t \in [0, 1], \label{4-12}\\
    &u(x, 0) = 3sech^{2}(0.5x),  \quad x \in [-40, 40] .\label{4-13}
\end{align}
Problems\eqref{4-11}-\eqref{4-13} have the exact solution  $u(x, t) =3sech^{2}(0.5x-0.5t)$.
\end{example}

For DDC-PINNs, the spatial domain $[-40,40]$ is decomposed into $[-40,0]$ and $[0,40]$ , with loss weights of $\lambda_{1}=1$, $\lambda_{2}=10$, $\lambda_{3}=1$, $\lambda_{4}=1$, $\lambda_{5}=10$, $\lambda_{6}=0.001$, $\lambda_{7}=0.01$, $\lambda_{8}=0.01$, and $\lambda_{9}=0.001$.
The setup of training point allocation for XPINNs, IDPINNs, and DDC-PINNs (first stage) is as follows:
\begin{itemize}
    \item A total of 3000 training points are assigned to each subdomain.
    \item An additional 1000 points are placed along the subdomain interface.
    \item Boundary conditions are enforced using 200 points.
    \item Initial conditions are imposed via 200 points.
\end{itemize}
For the CEENs method, the neural Galerkin method, and the second  stage of DDC-PINNs, the spatial discretization employs 200 points, and the time step size is set to 0.001.

Figure~\ref{fig:6} presents the exact solution, the DDC-PINNs prediction, and the corresponding two-dimensional and three-dimensional visualizations of the point-wise errors for problem \eqref{4-13}-\eqref{4-15}. The solution maintains a localized solitary-wave structure during propagation.
Figure  \ref{fig:7} compares the DDC-PINNs solution (blue solid line) with the exact solution (red dashed line) at different time instances.

In Table \ref{tab:5}, we summarize the relative $L^2$ errors and runtimes obtained by PINNs and the domain-decomposition methods (XPINNs, IDPINNs, and DDC-PINNs) for Example~\ref{ex:3}. Among the compared methods, DDC-PINNs produce the smallest mean relative $L^2$ error. In terms of runtime, DDC-PINNs require slightly more computational time than IDPINNs and substantially more computational time than PINNs and XPINNs in this example.
In Table \ref{tab:6}, we compare the relative $L^2$ errors and runtimes of the Neural Galerkin method, the CEENs method, and the DDC-PINNs method. For this example, DDC-PINNs achieve a lower mean relative $L^2$ error than both Neural Galerkin and CEENs. Regarding runtime, the Neural Galerkin method is faster than DDC-PINNs, whereas CEENs require a considerably longer computational time.

\begin{table}[h]
\centering
\caption{Relative $L^2$ error comparison among PINNs, XPINNs, IDPINNs, and DDC-PINNs for Example~\ref{ex:3}.}
 \label{tab:5}
\begin{tabular}{ccccc}
\toprule
\textbf{Methods} &PINNs & XPINNs & IDPINNs  & DDC-PINNs \\
\midrule
Relative \(L^2\) error  & \((3.15 \pm 1.12)\%\) & \((1.50\pm 0.45)\%\)  &\((0.33 \pm 0.16)\%\) &   \((0.20 \pm 0.08)\%\) \\
Runtime (s)
&
$391.18\pm 7.65$
&
$662.96\pm 7.63$
&
$913.34\pm 14.30$
&
$953.14\pm 21.99 $
\\
\bottomrule
\end{tabular}
\end{table}

\begin{table}[h]
\centering
\caption{Accuracy and runtime comparison of Neural Galerkin, CEENs, and DDC-PINNs for Example~\ref{ex:3}.}
\label{tab:6}
\begin{tabular}{cccc}
\toprule
Methods
&
Neural Galerkin
&
CEENs
&
DDC-PINNs
\\
\midrule
Relative $L^2$ error
&
$ (0.43 \pm 0.16)\%$
&
$ (13.52\pm 0.16)\%$
&
$  (0.20 \pm 0.08)\%$
\\
Runtime (s)
&
$171.70\pm 6.54$
&
$6587.25\pm 113.41$
&
$953.14\pm 21.99$
\\
\bottomrule
\end{tabular}
\end{table}

\begin{figure}[H]
    \centering
    \begin{subfigure}[b]{0.4\textwidth}
        \centering
        \includegraphics[width=\textwidth]{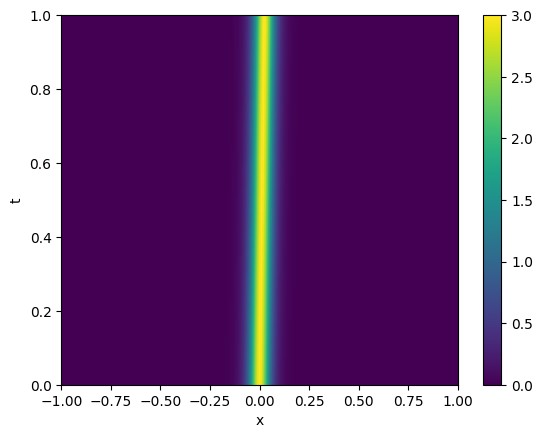}
        \caption{Exact solution}
    \end{subfigure}
    \vspace{0.5cm}
    \begin{subfigure}[b]{0.4\textwidth}
        \centering
        \includegraphics[width=\textwidth]{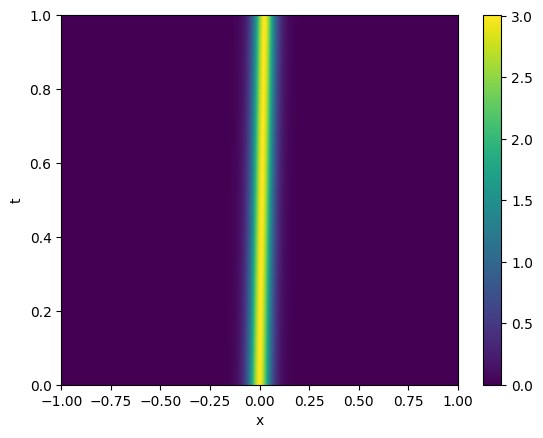}
        \caption{DDC-PINNs solution}
    \end{subfigure}
    \vspace{0.5cm}
    \begin{subfigure}[b]{0.4\textwidth}
        \centering
        \includegraphics[width=\textwidth]{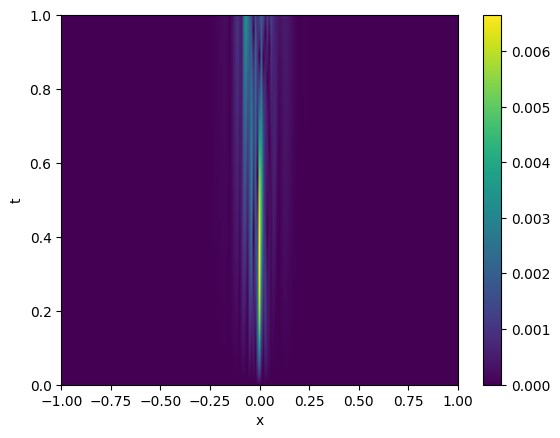}
        \caption{Point-wise errors}
    \end{subfigure}
    \vspace{0.5cm}
    \begin{subfigure}[b]{0.4\textwidth}
        \centering
        \includegraphics[width=\textwidth]{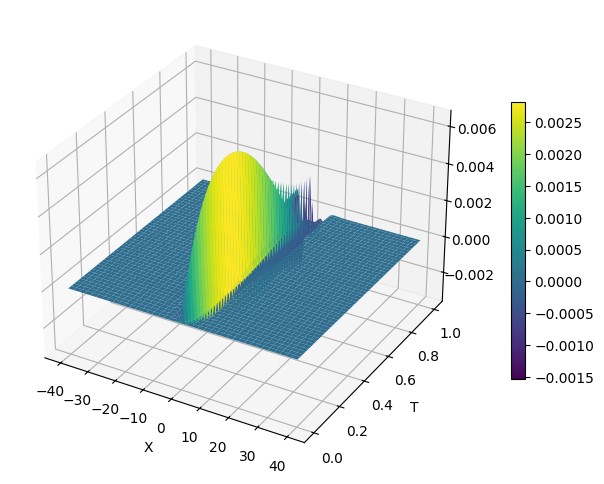}
        \caption{3D point-wise errors}
        \end{subfigure}
    \caption{Comparison of the exact solution, the DDC-PINNs prediction, and the corresponding point-wise errors for Example~\ref{ex:3}. Panels (c) and (d) present the two-dimensional point-wise absolute errors and the three-dimensional point-wise errors, respectively. The solution maintains a localized solitary-wave structure during propagation.}\label{fig:6}
\end{figure}

\begin{figure}[H]
    \centering
    \begin{subfigure}[b]{1\textwidth}
        \centering
        \includegraphics[width=\textwidth]{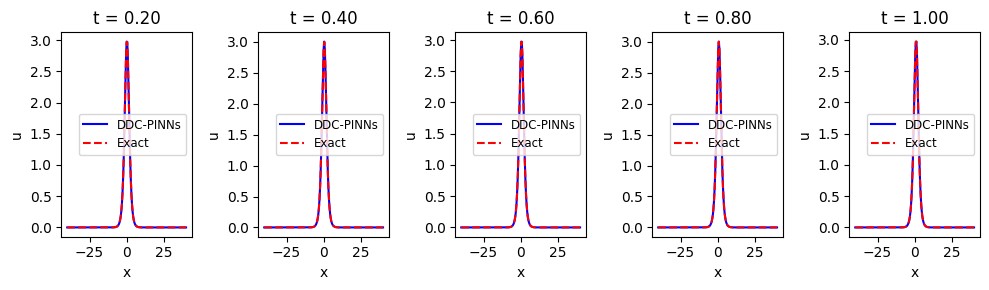}
    \end{subfigure}
    \caption{ Comparison between the exact solution (red dashed line) and the DDC-PINNs solution (blue solid line) at representative time instants for Example \ref{ex:3}.
The results demonstrate that both the temporal evolution and the spatial profile of the solitary wave are accurately captured by the proposed framework.}      \label{fig:7}
\end{figure}

\subsection{Two-Dimensional Heat Equation}
\begin{example}\label{ex:4}
This example considers a two-dimensional heat-conduction problem.
\begin{align}
    &u_t=u_{xx}+u_{yy}
     +e^{-t}(\sin( x)-\sin(y)) , \quad (x,y) \in \Omega, t \in [0, 1], \label{4-14}\\
    &u(x,y, t) =  0,  \quad (x,y) \in \partial \Omega   \quad t \in [0, 1], \label{4-15}\\
    &u(x,y, 0) = \sin(x) \sin(y),  \quad (x,y) \in \Omega .\label{4-16}
\end{align}
$\Omega= [0, \pi] \times [0, \pi]$,
The exact solution to problem \eqref{4-14}--\eqref{4-16} takes the form $u(x,y, t) = e^{-t}\sin(x)\sin(y)$.
\end{example}
The loss function weights used in DDC-PINNs are $\lambda_{1}=10$, $\lambda_{2}=5$, $\lambda_{3}=5$, $\lambda_{4}=0.03$, $\lambda_{5}=5$, $\lambda_{6}=0.03$, $\lambda_{7}=5$, $\lambda_{8}=5$, and $\lambda_{9}=0.03$. The spatial domain $[0,\pi] \times [0,\pi]$ is divided into two subdomains, $[0,\pi/2] \times [0,\pi]$ and $[\pi/2,\pi] \times [0,\pi]$.
For the three methods under consideration—XPINNs, IDPINNs, and DDC-PINNs (first stage)—the training point allocation is configured as:
\begin{itemize}
    \item Each subdomain: 5000 points.
    \item Subdomain interface: 1000 points.
    \item Boundary conditions: 6400 points.
    \item Initial conditions: 1600 points.
\end{itemize}
For the CEENs method, the neural Galerkin method, and the second  stage of DDC-PINNs, the spatial discretization uses 2500 points, and the time step is fixed at 0.001.

As shown in Figure~\ref{fig:8}, the figure presents the exact solution, the DDC-PINNs solution, and the two-dimensional and three-dimensional point-wise errors at $t=0.4$ and $t=0.8$ for problems~\eqref{4-14}-\eqref{4-16}.

In Table~\ref{tab:7}, we summarize the relative $L^2$ errors and runtimes obtained by PINNs, XPINNs, IDPINNs, and DDC-PINNs for Example~\ref{ex:4}. Among the four methods, DDC-PINNs produce the smallest mean relative $L^2$ error. In terms of runtime, DDC-PINNs require more computational time than PINNs, XPINNs, and IDPINNs for this example. The reported standard deviations indicate that all methods exhibit some variability across repeated runs.
In Table~\ref{tab:8}, we compare the relative $L^2$ errors and runtimes of the Neural Galerkin method, the CEENs method, and the DDC-PINNs method. For this example, the Neural Galerkin method achieves the smallest mean relative $L^2$ error, while DDC-PINNs obtain a lower mean relative $L^2$ error than CEENs. Regarding runtime, the Neural Galerkin method is slightly faster than DDC-PINNs, whereas CEENs require substantially longer computational time.


\begin{table}[h]
\centering
\caption{Relative $L^2$ error comparison among PINNs, XPINNs, IDPINNs, and DDC-PINNs for Example~\ref{ex:4}.}
 \label{tab:7}
\begin{tabular}{ccccc}
\toprule
\textbf{Methods} & PINNs & XPINNs & IDPINNs  & DDC-PINNs \\
\midrule
Relative \(L^2\) error  &  \((4.32 \pm 0.70)\%\) & \((7.28 \pm 1.06)\%\) &  \((1.32 \pm 1.02)\%\) &  \((0.72 \pm 0.35)\%\) \\
Runtime (s)
&
$1448.46\pm 19.49$
&
$1714.05\pm 18.52$
&
$2468.28\pm 36.76$
&
$2570.93\pm 70.03$
\\
\bottomrule
\end{tabular}
\end{table}

\begin{table}[h]
\centering
\caption{Accuracy and runtime comparison of Neural Galerkin, CEENs, and DDC-PINNs for Example~\ref{ex:4}.}
\label{tab:8}
\begin{tabular}{cccc}
\toprule
Methods
&
Neural Galerkin
&
CEENs
&
DDC-PINNs
\\
\midrule
Relative $L^2$ error
&
$(0.11\pm 0.02)\%$
&
$(50.19 \pm 0.25)\%$
&
$(0.72 \pm 0.35)\%$
\\
Runtime (s)
&
$2427.39\pm 41.94$
&
$16010.14 \pm 566.97$
&
$2570.93\pm 70.03$
\\
\bottomrule
\end{tabular}
\end{table}

\begin{figure}[h]
    \centering
    \begin{subfigure}[b]{0.35\textwidth}
        \centering
        \includegraphics[width=\textwidth]{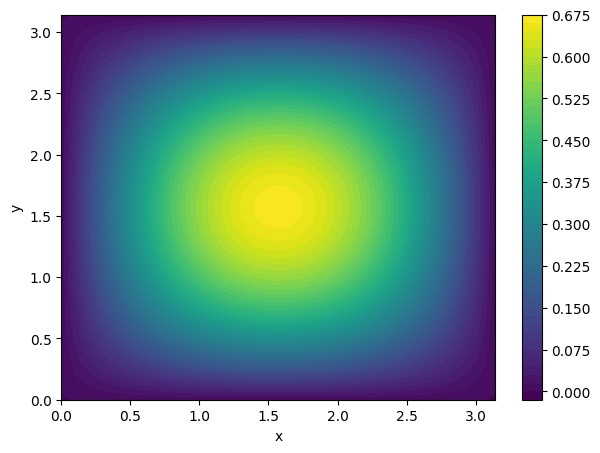}
        \caption{DDC-PINNs solution at $t=0.4$}
    \end{subfigure}
    \vspace{0.5cm}
    \begin{subfigure}[b]{0.35\textwidth}
        \centering
        \includegraphics[width=\textwidth]{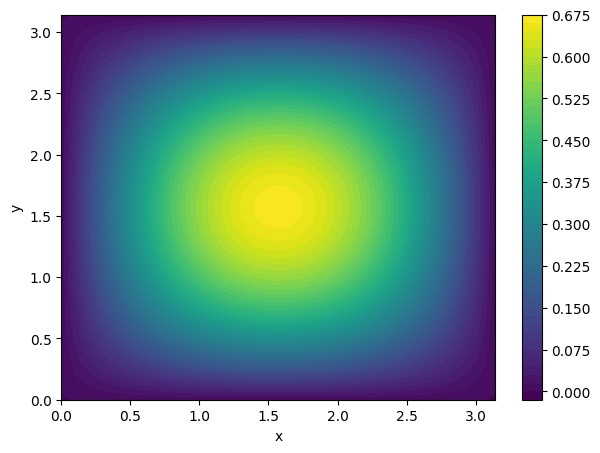}
        \caption{Exact solution at $t=0.4$}
    \end{subfigure}
    \vspace{0.5cm}
    \begin{subfigure}[b]{0.35\textwidth}
        \centering
        \includegraphics[width=\textwidth]{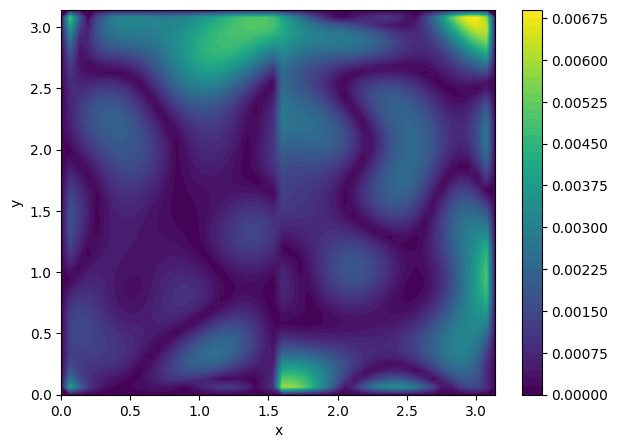}
        \caption{Point-wise errors at $t=0.4$}
    \end{subfigure}
        \vspace{0.5cm}
    \begin{subfigure}[b]{0.35\textwidth}
        \centering
        \includegraphics[width=\textwidth]{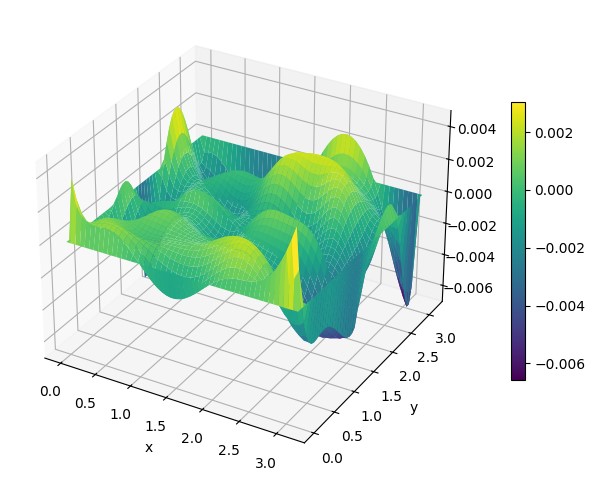}
        \caption{3D point-wise errors at $t=0.4$}
    \end{subfigure}
    
\vspace{0.5cm}
    \begin{subfigure}[b]{0.35\textwidth}
        \centering
        \includegraphics[width=\textwidth]{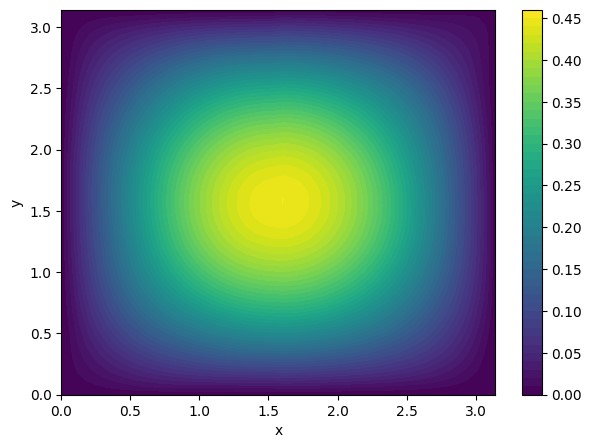}
        \caption{DDC-PINNs solution at $t=0.8$}
    \end{subfigure}
    \vspace{0.5cm}
    \begin{subfigure}[b]{0.35\textwidth}
        \centering
        \includegraphics[width=\textwidth]{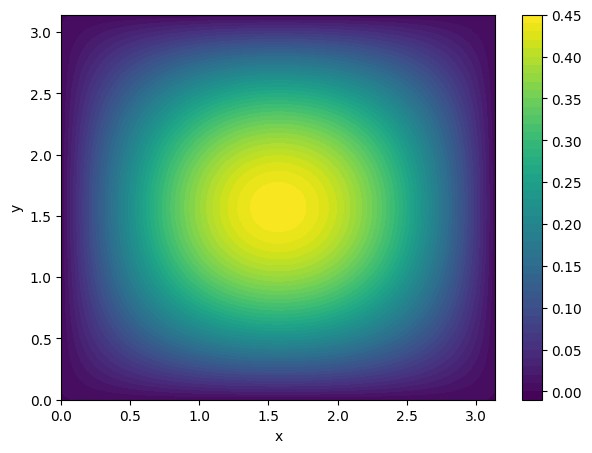}
        \caption{Exact solution at $t=0.8$}
    \end{subfigure}
    \vspace{0.5cm}
    \begin{subfigure}[b]{0.35\textwidth}
        \centering
        \includegraphics[width=\textwidth]{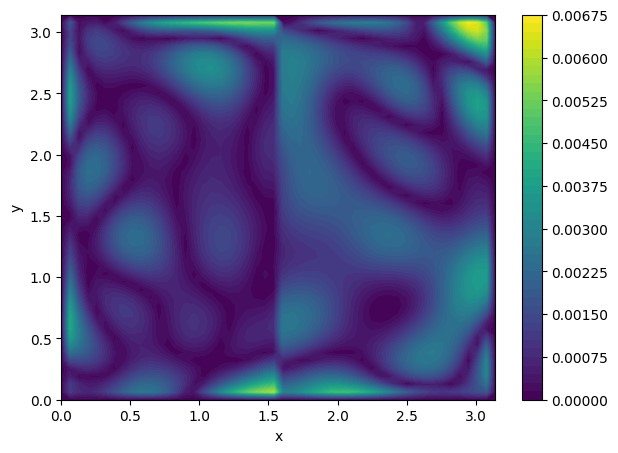}
        \caption{Point-wise errors at $t=0.8$}
    \end{subfigure}  
        \vspace{0.5cm}
    \begin{subfigure}[b]{0.35\textwidth}
        \centering
        \includegraphics[width=\textwidth]{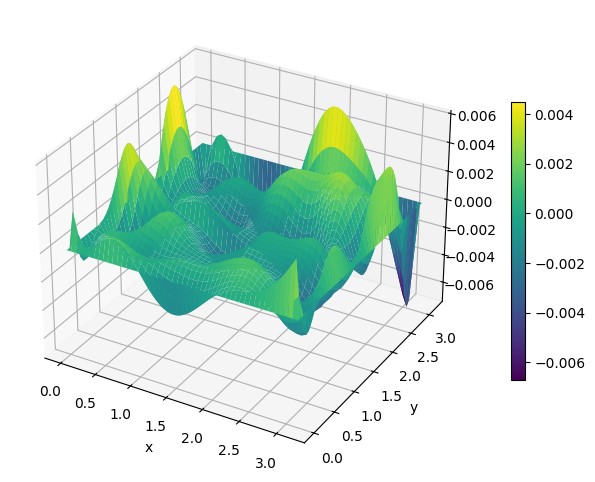}
        \caption{3D point-wise errors at $t=0.8$}
    \end{subfigure}
    \caption{Presents the exact solution, the DDC-PINNs prediction, and the corresponding two-dimensional and three-dimensional visualizations of the point-wise errors for Example~\ref{ex:4} at $t=0.4$ and $t=0.8$. The solution exhibits the expected two-dimensional diffusion behavior, with the amplitude decreasing over time.} \label{fig:8}
\end{figure}

\subsection{Hyperparameter Sensitivity Analysis}

To further investigate the influence of training hyperparameters on the proposed DDC-PINNs framework, we conduct an additional sensitivity analysis. Among the benchmark problems considered in this work, Example~\ref{ex:3} is selected as a representative test case because it involves nonlinear dynamics and higher-order spatial derivatives, making it more challenging than the diffusion-type examples.

In this study, the network architecture, loss-function weights, collocation points, initialization strategy, and all other settings are kept unchanged. Only the learning rate and the number of training iterations are varied. The purpose of this experiment is to assess the influence of these training hyperparameters on the prediction accuracy and computational cost of DDC-PINNs.

\begin{table}[htbp]
\centering
\caption{Sensitivity of DDC-PINNs to learning rate and training iterations for Example~\ref{ex:3}.}
\label{tab10}
\begin{tabular}{cccc}
\hline
Learning rate & Iterations & Relative $L^2$ error & Runtime (s) \\
\hline
$5\times10^{-4}$ & 20000 & $9.00\times10^{-4}$& 902.28 \\
$1\times10^{-3}$ & 20000 &$1.33 \times10^{-3}$ & 905.95\\
$2\times10^{-3}$ & 20000 &$1.12 \times10^{-3}$& 890.30\\
$4\times10^{-3}$ & 20000 &$1.49 \times10^{-3}$& 898.15\\
\hline
$1\times10^{-3}$ & 5000 & $4.92\times10^{-3}$& 238.61\\
$1\times10^{-3}$ & 10000 & $1.68\times10^{-3}$& 454.38 \\
$1\times10^{-3}$ & 20000 & $1.33\times10^{-3}$& 921.48\\
$1\times10^{-3}$ & 30000 &$2.09\times10^{-3}$& 1431.81 \\
$1\times10^{-3}$ & 40000 &$1.07 \times10^{-3}$ & 1843.97 \\
\hline
\end{tabular}
\end{table}

The results show that both the learning rate and the number of training iterations influence the final prediction accuracy and runtime. For the learning rates considered in this study, the relative $L^2$ errors remain of the same order of magnitude, indicating that the proposed framework is not highly sensitive to moderate changes in the learning rate. The smallest error is obtained with a learning rate of $5\times10^{-4}$, whereas larger learning rates lead to comparable accuracy with slight variations in the final error.

For a fixed learning rate of $10^{-3}$, increasing the number of training iterations generally reduces the prediction error from 5000 to 20000 iterations. Further increasing the number of iterations does not lead to a monotonic improvement in accuracy, although the smallest error among the tested settings is obtained at 40000 iterations. As expected, the runtime increases with the number of training iterations.

Overall, the results suggest that DDC-PINNs maintains comparable prediction accuracy over a reasonable range of hyperparameter choices. The hyperparameter setting adopted in the main experiments (learning rate $=10^{-3}$ and 20000 training iterations) provides a reasonable compromise between prediction accuracy and computational cost for the benchmark problems considered in this work.

\section{Conclusions}\label{sec6}

This paper proposes DDC-PINNs, a two-stage framework that combines domain-decomposition PINNs with classical time integration for solving time-dependent PDEs. By introducing additional smoothing constraints on both the solution and PDE residual, the proposed method improves the continuity of the neural approximation across subdomain interfaces. Numerical experiments on the diffusion equation, inviscid Burgers equation, KdV equation, and two-dimensional heat equation show that DDC-PINNs achieves competitive accuracy on the benchmark problems considered in this work. In three of the four examples, DDC-PINNs attains the smallest mean relative $L^2$ error among the compared methods, while in the two-dimensional heat equation example, the Neural Galerkin method yields a lower mean relative $L^2$ error. These results indicate that the proposed framework provides a viable approach for combining domain-decomposition neural approximations with sequential temporal evolution.
The main advantages of DDC-PINNs are summarized as follows:
(1) Domain decomposition with enhanced interface regularization improves the continuity and accuracy of the neural approximation across subdomain interfaces.
(2) The operator-freezing strategy reformulates the original PDE into an auxiliary ODE system, thereby separating spatial approximation from temporal evolution.
(3) Temporal evolution is performed without repeated neural-network optimization, preserving the temporal evolution order through sequential ODE integration.
The current limitations of DDC-PINNs include:
(1) The numerical experiments are limited to benchmark problems with relatively smooth solutions and prescribed domain decompositions.
(2) The applicability of the framework to complex geometries, higher-dimensional problems, stiff and multi-scale systems, and discontinuous solutions requires further investigation.
Therefore, the present study should be regarded as a proof-of-concept validation of the proposed framework. Future work will focus on adaptive domain decomposition strategies, more challenging PDEs, inverse problems, noisy observational data, and further theoretical analysis.

\section*{CRediT authorship contribution statement}
\textbf{Xun Yang}: Writing-review \& editing, Writing-original draft, Methodology, Investigation.
{\color{mgq} \textbf{Guanqiu Ma}: Writing-review \& editing, Supervision, Methodology.}
\textbf{Maohua Ran}: Writing-review \& editing, Supervision, Methodology, Funding acquisition.

\section*{Declaration of competing interest}
 The authors declare that they have no known competing financial interests or personal relationships that could have appeared to influence the work reported in this paper.

\section*{Acknowledgments}\label{sec7}
This work is supported by the Sichuan Science and Technology Programs (No. 2024NSFSC0441).

\section*{Data availability statement}
 All source code for this manuscript is openly shared on Zenodo at  \url{https://doi.org/10.5281/zenodo.20628188}


\end{document}